\documentclass[12pt]{article}
\usepackage{amsmath}
\usepackage{amssymb}
\newcommand{\Z}{\mathbb Z}

\newcommand{\N}{\mathbb N}
\newcommand{\C}{\mathbb C}

\newcommand{\p}{\partial}

\newcommand{\La}{\Lambda}
\newcommand{\la}{\lambda}
\newcommand{\f}{\phi}
\newcommand{\bV}{\overline{V}}

\begin{document}

\title{Mathieu-Zhao spaces of polynomial rings}
\author{Arno van den Essen$ $\footnote{corresponding author: essen@math.ru.nl}\,\, and Loes van Hove}
\maketitle

\begin{abstract} We describe all Mathieu-Zhao spaces of $k[x_1,\cdots,x_n]$ ($k$ is an
algebraically closed field of characteristic zero) which contains an ideal of finite codimension. Furthermore we give
an algorithm to decide if a subspace of the form $I+kv_1+\cdots+kv_r$ is a Mathieu-Zhao space, in case
the ideal $I$ has finite codimension.
\end{abstract}

\section*{Introduction}

Since its formulation in 1939 by Keller the Jacobian Conjecture has been studied by many authors, but 
remains open in all dimensions greater than one. Many attempts have been made to generalize this conjecture,
however most of these generalizations turned out to be false. Only one such a conjecture, due to Olivier Mathieu 
in [6], is still open. More recently Wenhua Zhao came up with several amazing new conjectures, all implying
the Jacobian Conjecture. Even better, he created a new framework in which all these fascinating conjectures, including
Mathieu's conjecture, can be studied:this is his theory of Mathieu subspaces ([7], [8 ], [9 ], [10] and [1]). The name Mathieu subspaces was recently changed into Mathieu-Zhao spaces, for short MZ-spaces, by the first author in [2].

An MZ-space is a  generalization of the notion of an ideal in a ring. More precisely, let $k$ be a field,
$R$ a $k$-algebra and $V$  a $k$-linear subspace of $R$. Then $V$ is called a {\em (left) MZ-space
of $R$} if the following holds: if $a\in R$ is such that $a^m\in V$, for all large $m$ (i.e.\@ there exists $N$ such that
$a^m\in V$ for all $m\geq N$), then for all $b\in R$ also $ba^m\in V$ for all large $m$.

The new conjectures  introduced by Zhao all concern MZ-spaces of polynomial rings over a field.
Therefore one is naturally led  to the  study of MZ-spaces of such rings. A first step toward 
a description of these spaces, for the case of univariate polynomial rings, was made in [3]. There the authors
classify all MZ-spaces of $k[t]$ which contain a non-zero ideal. These spaces have finite codimension. However
classifying MZ-spaces, even of codimension one of $k[t]$, is still far too complicated. For example the set of all $f\in\C[t]$
such that $\int_0^1 f(t)\,dt=0$ is an MZ-space of $\C[t]$, but its proof is not at all obvious  (see for example [4] or [1]).

The aim of this paper is to extend the results obtained in [3] to polynomial rings in $n$ variables. More precisely, in case
$k$ is an algebraically closed field of characteristic zero,
we give a complete description of all MZ-spaces of $k[x]:=k[x_1,\cdots,x_n]$ containing an ideal of finite
codimension. Furthermore, we give an algorithm which decides if a given subspace
of $k[x]$ of the form $I+kv_1+\cdots+kv_h$ is an MZ-space, in case $I$ has finite codimension. 

The results described in this paper where first obtained by the second author in her Master's thesis [5], at the Radboud University in Nijmegen. This
paper contains some simplifications of the original proofs.

\section{Preliminaries and notations}

Throughout this paper $k$ will denote an algebraically closed field of characteristic zero and $k[x]:=k[x_1,\cdots,x_n]$
is the polynomial ring in $n$ variables over $k$. $V$ will always denote a $k$-linear subspace of $k[x]$ and we
additionally {\em assume} that $V$ contains an ideal $I$ such that $k[x]/I$ is a finite dimensional $k$-vectorspace, say
of dimension $d$. It follows that the $d+1$ vectors $1,\overline{x_1},\overline{x_1^2},\cdots, \overline{x_1^d}$
are linearly dependent over $k$, which implies that $I$ contains a monic polynomial $f_1(x_1)\in k[x_1]$ of degree say $d_1\geq 1$.
Since this argument can be repeated for every $i$ we deduce that there exist monic polynomials $f_1(x_1),\cdots,f_n(x_n)$, 
of positive degrees $d_1,\cdots,d_n$ respectively, such that $J:=(f_1(x_1),\cdots,f_n(x_n))\subseteq I\subseteq V$.
Observe that dim$_k k[x]/J=d_1\cdots d_n$ is finite. Consequently we may, and will assume from now on that $I=(f_1(x_1),\cdots,f_n(x_n))$. 

The advantage of this assumption is that $A:=k[x]/I$ has a nice structure. To see this let's fix some notations. First we 
denote by $\La_i$ the set of different zeros of $f_i$ in $k$ and for $\la_i\in\La_i$ we denote by $m(\la_i)$ its
multiplicity. So
$$f_i(x_i)=\prod_{\la_i\in\La_i}(x_i-\la_i)^{m(\la_i)}$$
\noindent We may assume that $0\notin\La_i$ for all $i$: just replace $x_i$ by $x_i-c_i$ for some suitable $c_i\in k$ and
observe that sending each $x_i$ to $x_i-c_i$ is a $k$-automorphism of $k[x]$.
\noindent Now define $\La=\La_1\times\cdots\times\La_n$. So an element $\la\in\La$ is an $n$-tuple
of the form $\la=(\la_1,\cdots,\la_n)$, where each $\la_i$ belongs to $\La_i$. The $n$-tuple $(m(\la_1),\cdots,m(\la_n))$ we denote by $m(\la)$. If furthermore for each $\la\in\La$ we
denote by $[(x-\la)]^{m(\la)}$ the ideal $((x_1-\la_1)^{m(\la_1)},\cdots,(x_n-\la_n)^{m(\la_n)})$ in $k[x]$, it follows from the Chinese remainder theorem and an easy induction that
$$k[x]/I\simeq\prod_{\la\in\La} k[x]/[(x-\la)]^{m(\la)}$$
The isomorphism is given by $\f(g+I)=(g+[(x-\la)]^{m(\la)})_{\la\in\La}$. The ring on the right-hand side
we denote by $B$. It is a product of the local rings $B_{\la}:= k[x]/[(x-\la)]^{m(\la)}$. Hence each such a ring has only
two idempotents, namely $0$ and $1$. It follows
that the elements $e_{\la}=(0,\cdots,0,1,0,\cdots,0)\in B$ (where the $1$ appears at  the component with index $\la$)
form an {\em orthogonal basis of idempotents of $B$}, i.e.\@ each $e_{\la}$ is a non-zero idempotent of $B$, 
$e_{\la}\cdot e_{\mu}=0$ for all $\la\neq\mu\in\La$ and each non-zero idempotent of $B$ is of the form $\sum_{\la\in\La^{'}}e_{\la}$,
for some non-empty subset $\La^{'}$ of $\La$. By the isomorphism $\f$ there exist $g_{\la}\in k[x]$, such that
$\f(g_{\la}+I)=e_{\la}$. Consequently the elements $g_{\la}+I$ form an orthogonal bases of idempotents of $A$.

To understand the importance of these idempotents we recall two facts from [10]. The first fact says  that $V$ is an MZ-space of $k[x]$ if and only if $\bV:=V/I$ is an MZ-space
of $A$. So we need to study MZ-spaces of $A$. Therefore observe that $A$ is finite dimensional over $k$, so all its elements are algebraic over $k$. It then follows
from Zhao's idempotency theorem (theorem 4.2, [10]) that $\bV$ is an MZ-space of $A$ if and only if for each idempotent
$e$ of $A$, which belongs to $\bV$, the ideal $Ae$ is contained in $\bV$. Before we can use these results to obtain
a first characterization of MZ-spaces of $k[x]$ containing $I$, we need one more result,
which will be applied to the ring $A$ and the idempotents $g_{\la}+I$ described above:

\medskip

\noindent{\bf Lemma 1.} {\em Let $R$ be a commutative ring which has an orthogonal
basis $E$ of idempotents.  If $M$ is an MZ-space of $R$, 
 then the only idempotents of $R$ in $M$ are $0$ or the elements of the form $\sum_j e_j$, where each $e_j$ belongs to $E_0:=E\cap M$.}

\medskip

\noindent{\em Proof.} Let $e\in M$ be an idempotent and assume that $e\neq 0$. Then $e=\sum_j e_j$, for some $e_j\in E$. Assume that one of these $e_j$ does not belong to $E_0$, say $e_i\notin E_0$. Then $e_i\notin M$.  Now observe 
that $e^m=e\in M$ for all $m>0$. Since $M$ is an MZ-space this implies that $e_ie^m\in M$ for all large $m$.
However $e_ie^m=e_ie=e_i^2=e_i$. So $e_i\in M$, a contradiction. So  $e_j\in E_0$, for each $j$.

\medskip

Now we are able to prove the first main theorem. Therefore let $\La_0$ be the set of $\la\in\La$ such that
$g_{\la}\in V$. Furhermore, for each $\La'\subseteq\La$ we put  $I(\La')=\sum_{\la\in\La'}g_{\la}$ if $\La'\neq \emptyset$ and $I(\La')=0$ otherwise.

\medskip

\noindent{\bf Theorem 1.} {\em $V$ is an MZ-space of $k[x]$ if and only if for each non-empty subset $\La'$ of $\La$
the following conditions hold:\\
\noindent i) $I(\La'\backslash\La_0)\notin V$, if $\La\backslash\La_0\neq\emptyset$.\\
\noindent ii) $k[x]\cdot I(\La'\cap\La_0)\subseteq V$.}

\medskip

\noindent{\em Proof.} $(\Rightarrow)$ Assume  $\La'\backslash\La_0\neq\emptyset$. Suppose that 
$\sum_{\la\in\La'\backslash\La_0}g_{\la}\in V$. Then $\sum_{\la\in\La'\backslash\La_0}\overline{g_{\la}}\in \bV$. Since
$V$ is an MZ-space of $k[x]$, $\bV$ is an MZ-space in $A$. So by lemma 1 (applied to the ring $A$ and the idempotents
$g_{\la}+I$) it follows that $\overline{g_{\la}}=g_{\la}+I\in\bV$, for all $\la\in\La'\backslash\La_0$. Since $I\subseteq V$, this implies that
$g_{\la}\in V$ for all these $\la$. However if $\la\in\La'\backslash\La_0$, then in particular $\la\notin\La_0$.
So $g_{\la}\notin V$, contradiction. This proves i). To see ii) just observe that  $\overline{I(\La'\cap\La_0)}=
\sum_{\la\in\La'\cap\La_0}\overline{g_{\la}}$ is an idempotent in $A$ which is contained in $\bV$. Since $\bV$ is an MZ-space in $A$ (for $V$ is one in $k[x]$),
it follows from Zhao's idempotency theorem that $A\cdot \overline{I(\La'\cap\La_0)}\subseteq\bV$. Using again that
$I\subseteq V$ this implies ii).\\
\noindent $(\Leftarrow)$ It suffices to show that $\bV$ is an MZ-space of $A$. We use Zhao's idempotency theorem.
So let $e\in\bV$ be a non-zero idempotent of $A$. Then there exists a non-empty subset $\La'$ of $\La$ such that
$e=\sum_{\la\in\La'}\overline{g_{\la}}\in \bV$. Split this sum into
$$\sum_{\la\in\La'\backslash\La_0}\overline{g_{\la}}+\sum_{\la\in\La'\cap\La_0}\overline{g_{\la}}$$
\noindent By definition of $\La_0$ the last part belongs to $\bV$. Consequently $\sum_{\la\in\La'\backslash\La_0}\overline{g_{\la}}\in\bV$, whence $\sum_{\la\in\La'\backslash\La_0}g_{\la}\in V$.
It follows from i) that $\La'\backslash\La_0=\emptyset$. So each non-zero idempotent of $\bV$ is  of the form
$\sum_{\la\in\La'\cap\La_0}\overline{g_{\la}}$. By ii) we get that $A\cdot\sum_{\la\in\La'\cap\La_0}\overline{g_{\la}}\subseteq\bV$. So by Zhao's idempotency theorem we deduce
that $\bV$ is an MZ-space of $A$, which completes the proof.

\section{$V$ as the kernel of a linear map}

We recall that $V$ is a $k$-linear subspace of $k[x]$ containing an ideal $I$ of the form $I=(f_1(x_1),\cdots,f_n(x_n))$,
where each $f_i$ is a univariate polynomial of positive degree $d_i$. It follows that $A:=k[x]/I$ is
finite dimensional over $k$ and hence so is $k[x]/V$. If $r$ denotes the dimension of this space, there exists a $k$-linear
isomorphism $\psi:k[x]/V\rightarrow k^r$. Let $\pi$ be the canonical map 
from $k[x]$ to $k[x]/V$. Then $\frak{L}:=\psi\circ\pi$ is a surjective $k$-linear map from $k[x]$ to $k^r$ such that
$V=ker\,\frak{L}$. Write $\frak{L}=(L_1,\cdots,L_r)$. Then each $L_i:k[x]\rightarrow k$ is a $k$-linear map
having $I$ in its kernel.
In the remainder of this section we give an explicit description  of such $k$-linear maps. In order to do so we introduce some
more notation: if $f\in k[x]$ we let
$$Deg\,f:=(deg_{x_1}f,\cdots,deg_{x_n}f)$$
\noindent and if $a,b\in\Z^n$ we define $a<b$ if and only if $a_i<b_i$ for all $i$. Furthermore we introduce two types of
operators on $k[x]$: the differential operators $D_j=x_j\p_{x_j}$, for each $j$ and the substitution
maps $S_{\la}:k[x]\rightarrow k$, given by $S_{\la}(g)=g(\la)$, for all $g\in k[x]$ and each $\la\in\La$. Finally
write $D:=(D_1,\cdots,D_n)$. With these notations we have:

\medskip

\noindent{\bf Theorem 2.} {\em Let $L:k[x]\rightarrow k$ be a $k$-linear map such that $I\subseteq ker\,L$. Then for every $\la\in\La$ there exists a polynomial $P_{\la}\in k[x]$ with $Deg\,P_{\la}<m(\la)$ such that
$L=\sum_{\la\in\La} S_{\la}\circ P_{\la}(D)$.}

\medskip

\noindent To prove this result we need some preparations:

\medskip

\noindent{\bf Lemma 2.} {\em $D_i^p(k[x](x_i-\la_i)^q)\subseteq k[x](x_i-\la_i)^{q-p}$, if $q>p\geq 0$.}

\medskip

\noindent{\em Proof.} Follows readily from Leibniz' rule and induction on $p$.

\medskip

\noindent{\bf Corollary.} {\em If $P_{\la}\in k[x]$ with $deg\,P_{\la}<m(\la)$, then 
$I\subseteq ker\,S_{\la}\circ P_{\la}(D)$.}

\medskip

\noindent{\em Proof.} We need to prove that $S_{\la}\circ P_{\la}(D)(k[x]f_i(x_i))=0$, for all $i$. We only treat the case
$i=1$. So let $a(x)\in k[x]$, we will show that $S_{\la}\circ P_{\la}(D)(a(x)f_1(x_1))=0$. Write $a(x)f_1(x_1)=b(x)(x_1-\la_1)^{m(\la_1)}$. Now observe that a typical monomial appearing in $P_{\la}(x)$ is of the form $cx_1^{i_1}\cdots x_n^{i_n}$, with $c\in k$
and $i_j<m(\la_j)$ for all $j$. So for the corresponding monomial in $P_{\la}(D)$ we get
$$cD_1^{i_1}\cdots D_n^{i_n}(a(x)f_1(x_1))=cD_2^{i_2}\cdots D_n^{i_n}D_1^{i_1}(b(x)(x_1-\la_1)^{m(\la_1)})$$
$$=_{lemma\,2}cD_2^{i_2}\cdots D_n^{i_n}(g(x)(x_1-\la_1)^{m(\la_1)-i_1}),\,\, g(x)\in k[x]$$
$$=cD_2^{i_2}\cdots D_n^{i_n}(g(x))(x_1-\la_1)^{m(\la_1)-i_1}$$
\noindent Since $i_1<m(\la_1)$ applying the substitution map $S_{\la}$ gives zero. Since this holds for every monomial
appearing in $P_{\la}(x)$, this completes the proof.

\medskip

\noindent{\em Proof of theorem 2.} If $L=0$, choose $P_{\la}=0$ for all $\la$. So let $L\neq 0$. Then there exists
$v\in k[x]$ with $L(v)=1$ and $k[x]/ker\,L\simeq k$. In particular $k[x]=ker\,L\oplus kv$. Since $I\subseteq ker\,L$
reduction modulo $I$ gives that $A=k[x]/I=\overline{ker\,L}\oplus k\overline{v}$. Let $d=dim_k A$. Choose
a $k$-basis $\overline{v_1},\cdots,\overline{v_{d-1}}$ of $\overline{ker\,L}$. Then $k[x]=I\oplus kv_1\oplus\cdots\oplus kv_{d-1}\oplus kv$.

For each $\la\in\La$ we define the universal polynomial
$$P_{\la}^U:=\sum_{i<m(\la)}P_{\la,i} x^i$$
\noindent where the $P_{\la,i}$ are variables. We will show that there exist $p_{\la,i}\in k$ such that $L$ equals
$L(p):=\sum_{\la\in\La} S_{\la}\circ (\sum_{i<m(\la)} p_{\la,i}D^i)$. Therefore we first observe that there are
$m(\la_1)\cdots m(\la_n)$ monomials $x^i$ with $i<m(\la)$. Hence there are $m(\la_1)\cdots m(\la_n)$ corresponding
variables $P_{\la,i}$. So summing over all $\la\in\La$ we get
$$\sum_{\la_1\in\La_1}\cdots\sum_{\la_n\in\La_n} m(\la_1)\cdots m(\la_n)=\sum_{\la_1\in\La_1}m(\la_1)\cdots\sum_{\la_n\in\La_n}m(\la_n)=d_1\cdots d_n$$
\noindent variables, which is precisely $d$, the dimension of $k[x]/I$. From the corollary above we know that for
each choice of the $p_{\la,i}\in k$ the corresponding operator $L(p)$ has $I$ in its kernel. Now we need to find  $p_{\la,i}\in k$ such that $L(p)$ is equal to $L$. Since the elements $v_1,\cdots,v_{d-1}$ belong to $ker\,L$
(for $\overline{v_i}\in \overline{ker\,L}$ and $I\subseteq ker\,L$), we must choose the $p_{\la,i}\in k$ in such a way
that $L(p)(v_i)=0$, for all $1\leq i\leq d-1$. This means that we have to solve a system of $d-1$ linear equations
in the $d$ variables $P_{\la,i}$. It follows that there exists at least one non-zero solution of $p_{\la,i}$'s in $k^d$. Let
$L(p)$ be the corresponding linear map. So $L$ and $L(p)$ are both zero on $I$ and the $v_i$. Since $k[x]=I\oplus kv_1\oplus\cdots\oplus kv_{d-1}\oplus kv$, it remains to see if they are equal on $v$. In general they are not. But we can
change the operator a little as follows:  define $a:=L(p)(v)$. We will
show below that $a\in k^*$. Since $L(v)=1$ it follows  that $L=(1/a)\cdot L(p)$ and $L$ not only agree on $I$
and the $v_i$ (where they both are zero), but also on $v$. So $L=(1/a)\cdot L(p)=L((1/a)p)$, which
completes the proof.

It remains to see that $a$ is non-zero. So assume that $a=0$. Then $L(p)$ is the zero-map, so $L(p)(x^m)=0$
for all monomials $x^m=x_1^{m_1}\cdots x_n^{m_n}$. From the definition of $L(p)$ and the fact that
$$D^i(x^m)=m_1^{i_1}\cdots m_n^{i_n}x^m$$
\noindent it then follows that
$$\sum_{(\la_1,\cdots,\la_n)\in\La}\sum_{i<m(\la)}p_{\la,i}m_1^{i_1}\cdots m_n^{i_n}\la_1^{m_1}\cdots\la_n^{m_n}=0$$
\noindent for all $(m_1,\cdots,m_n)\in\overline{\N}^n$. Then lemma 3 below gives that all  $p_{\la, i}$ are zero, a contradiction. So $a\neq 0$.

\medskip

\noindent{\bf Lemma 3.} {\em For each $i=(i_1,\cdots,i_n)\in\overline{\N}^n$ and $\la\in\La$ define
$f_{\la,i}:\overline{\N}\rightarrow k$ by
$$f_{\la,i}(m_1,\cdots,m_n)=m_1^{i_1}\cdots m_n^{i_n}\la_1^{m_1}\cdots\la_n^{m_n}$$
\noindent Then the $f_{\la,i}$ are linearly independent over $k$.}

\medskip

\noindent{\em Proof.} By induction on $n$. The case $n=1$ follows from the theory of linear recurrence relations (recall that
all $\la_i$ are non-zero).
So let $n\geq 2$ and assume that $\sum a_{\la,i}f_{\la,i}=0$, for some $a_{\la,i}\in k$. Then
$$\sum_{(i_n,\la_n)}(\sum_{(i',\la')} a_{\la,i}m_1^{i_1}\cdots m_{n-1}^{i_{n-1}}\la_1^{m_1}\cdots\la_{n-1}^{m_{n-1}})m_n^{i_n}\la_n^{m_n}=0$$
\noindent where $i'=(i_1,\cdots,i_{n-1})$ and $\la'=(\la_1,\cdots,\la_{n-1})$. From the case $n=1$ it then
follows that for each $i_n,\la_n$ the coefficent of the term $m_n^{i_n}\la_n^{m_n}$ equals zero, i.e.
$$\sum_{(i',\la')} a_{\la,i}m_1^{i_1}\cdots m_{n-1}^{i_{n-1}}\la_1^{m_1}\cdots\la_{n-1}^{m_{n-1}}=0$$
Then the induction hypothesis implies that all $a_{\la,i}$ are zero, which completes the proof.

\section{The main theorem}

Now we are able to give the main theorem of this paper. The notations are as introduced before. So $I=(f_1(x_1),\cdots,f_n(x_n))$ is contained in the $k$-linear subspace $V$ of $k[x]$ and the $g_{\la}+I$ form an orthogonal
basis of idempotents of $A=k[x]/I$. Furthermore $V=ker\,\frak{L}$, where $\frak{L}=(L_1,\cdots,L_r):k[x]\rightarrow k^r$ and each $L_i$ is of the form $L_i=\sum_{\la\in\La}S_{\la}\circ P_{\la}^{(i)}(D)$, for some
$P_{\la}^{(i)}\in k[x]$ with Deg $P_{\la}^{(i)}<m(\la)$, for all $\la$.

\medskip

\noindent{\bf Theorem 3.} {\em $V$ is an MZ-space of $k[x]$ if and only if the following two properties hold:\\
\noindent i) For each $\La'\subseteq\La$ such that 
$\La'\backslash\La_0\neq\emptyset$ there exists an $i$ such that
$$\sum_{\la\in\La'\backslash\La_0}P_{\la}^{(i)}(0)\neq 0$$
\noindent ii) $L_i(\sum_{\la\in\La'\cap\La_0}k[x]g_{\la})=0$, for all $1\leq i\leq r$.}

\medskip

\noindent{\em Proof.} By theorem 1 we know that $V$ is an MZ-space of $k[x]$ if and only if $I(\La'\backslash\La_0)\notin V$, when $\La'\backslash\La_0\neq\emptyset$ and $k[x]\cdot I(\La'\cap\La_0)\subseteq V$. The first condition
is equivalent to $\sum_{\la\in\La'\backslash\La_0}\frak{L}(g_{\la})\neq 0$, i.e.\@ to
$\sum_{\la\in\La'\backslash\La_0}L_i(g_{\la})\neq 0,\,\mbox{ for some } i$. By lemma 4 below $L_i(g_{\la})=P_{\la}^{(i)}(0)$,
which gives the first part of the theorem. The second condition $k[x]\cdot I(\La'\cap\La_0)\subseteq V$ is equivalent to
statement ii) of the theorem. This completes the proof.

\medskip

\noindent{\bf Lemma 4.} {\em Let $L$ be as in theorem 2. Then $L(g_{\la})=P_{\la}(0)$.}

\medskip

\noindent{\em Proof.} $L(g_{\la})=\sum_{\mu\in\La,\mu\neq\la}S_{\mu}\circ P_{\mu}(D)(g_{\la})+S_{\la}\circ P_{\la}(D)(g_{\la})$. Since by definition $g_{\la}\in [(x-\mu)]^{m(\mu)}$, for all $\mu\neq\la$, it follows from the fact that
Deg $P_{\mu}<m(\mu)$ that the first sum equals zero (copy the argument in the proof of the corollary above). So
$L(g_{\la})=S_{\la}\circ P_{\la}(D)(g_{\la})$. Finally, using the fact that $g_{\la}\equiv 1\,mod\,[(x-\la)]^{m(\la)}$ and
that $D_i(1)=0$ for all $i$, the result follows.

\section{Some final remarks}

\subsection*{An algorithm}

In the previous section we gave a complete description of the MZ-spaces of $k[x]$ containing an ideal of finite
codimension. It turned out that all these spaces are of the form
$$I+kv_1+\cdots +kv_h$$
\noindent where $I=(f_1(x_1),\cdots,f_n(x_n))$ and each $f_i(x_i)$ is an univariate polynomial
of positive degree.
As we will show now the results obtained above can also be used to give an algorithm which decides if a given space
of the form $I+kv_1+\cdots +kv_h$ is an MZ-space of $k[x]$,  when $I$ has finite codimension.

First, using Gr\"obner basis theory one can decide if  $I$ has finite codimension and in case it has find monic univariate
polynomials $f_i(x_i)$ of positive degrees $d_i$ contained in $I$.
As observed in the beginning of this paper, we can replace $I$ by the ideal generated by these $f_i(x_i)$. This also gives us the set $\La$. Next we need to determine the elements $g_{\la}$.
Since for each pair $\la,\mu\in\La$, with $\la\neq\mu$, the ideals $[(x-\la)]^{m(\la)}$ and $[(x-\mu)]^{m(\mu)}$ are
comaximal, we can find elements $i_{\la}\in [(x-\la)]^{m(\la)}$ and $i_{\mu}\in [(x-\mu)]^{m(\mu)}$ such that $i_{\la}+i_{\mu}=1$. Then one readily verifies that if we define
$$g_{\la}=\prod_{\mu\neq\la} i_{\mu}$$
\noindent  these elements have the desired properties.

Next we want to write $V$ as the kernel of a suitable linear map $\frak{L}$. Since the classes $\overline{x^m}$ with $m<(d_1,\cdots,d_n)$
form a basis of $k[x]/I$ it follows that the dimension of $k[x]/I$ equals $d:=d_1\cdots d_n$. Furthermore we can construct a $k$-basis of $\overline{V}:=V/I$. In other words replacing the original $v_i$
by better $v$'s we may assume that the elements $\overline{v_1}\cdots,\overline{v_h}$ form a $k$-basis of $\overline{V}$. Since $k[x]/I\big /V/I\simeq k[x]/V$ it follows that the dimension of $k[x]/V$ equals $r:=d-h$.

Then following the argument in the proof of theorem 2 one can construct a linear map $\frak{L}=(L_1,\cdots,L_r):k[x]\rightarrow k^r$, with ker $\frak{L}=V$ and each $L_i$ of the form as in theorem 2. Then to
decide if $V$ is an MZ-space of $k[x]$ we need to check the two properties given in theorem 3.

To do this we first compute $\La_0$, just by checking for which $\la\in\La$ we have $\frak{L}(g_{\la})=0$. The first condition of theorem 3  consist of a finite number of calculations, just one for each subset $\La'$ of $\La$ such that
$\La'\backslash\La_0\neq\emptyset$. Finally, the second condition $L_i(\sum_{\la\in\La'\cap\La_0}k[x]g_{\la})=0$, for all $1\leq i\leq r$, is equivalent to $L_i(\sum_{\la\in\La'\cap\La_0}x^mg_{\la})=0$, for all 
$1\leq i\leq r$  and all $m<(d_1,\cdots,d_n)$ (since each element of $k[x]$ is equivalent mod $I$ to a lineair combination
of monomials of the form $x^m$, with $m<(d_1,\cdots,d_n)$ and each $L_i$ has $I$ in its kernel). So again this only needs a finite number of calculations.

\subsection*{MZ-spaces of finitely generated Artin rings}

Let $R$ be a finitely generated $k$-algebra. Then $R$ is an Artin ring if and only if the dimension of $R$ is zero, or equivalently if $R$ is isomorphic to a quotient ring of the form $k[x_1,\cdots,x_n]/I$, for some $n$ and an ideal
$I$ of finite codimension. So studying MZ-spaces of $R$ amounts to studying MZ-spaces of $k[x]/I$, which in turn
amounts to studying MZ-spaces of $k[x]$ containing an ideal $I$ of finite codimension. This is exactly what
we did in the previous section. In other words, the main theorem of this paper completely describes
all MZ-spaces of Artin rings, which are finitely generated over $k$. Furthermore the algorithm given
above gives an algorithm to recognize MZ-spaces of R. 

\section*{References}

\noindent [1] A. van den Essen, {\em The Amazing Image Conjecture}, http://arxiv.org/abs/\\1006.5801 (2010).\\
\noindent [2] A. van den Essen, {\em An introduction to Mathieu subspaces}. Lectures delivered at the Chern Institute of Mathematics, Tianjin, China, July 2014.\\
\noindent [3] A. van den Essen and S. Nieman, {\em Mathieu-Zhao spaces of univariate polynomial rings with non-zero
strong radical}, J. Pure and Appl.
Algebra 220 (9) (2016), 3300-3306.\\
\noindent [4] J.P. Francoise, F. Pakovich, Y.Yomdin, W. Zhao, {\em Moment vanishing problem and positivity: Some
examples}, Bull. Sci. Math. 135 (2011), 10-32.\\
\noindent [5] L. van Hove, {\em Mathieu-Zhao subspaces}, Master's thesis University of Nijmegen, July 2015.\\
\noindent [6] O. Mathieu, {\em Some conjectures about invariant theory and their applications}, Alg\`ebra non commutative, groupes quantiques et invariants (Reims, 1995), S\'emin. Congr. 2, Soc. Math. France, 263-279 (1997).\\
\noindent [7] W. Zhao, {\em Hessian Nilpotent Polynomials and the Jacobian Conjecture}, Trans. Amer. Math. Soc. 359 (2007), 274-294.\\
\noindent [8] W. Zhao, {\em Images of commuting differential oparators of order one with constant leading coefficients},
J. Algebra 324 (2010), 231-247.\\
\noindent [9] W. Zhao, {\em Generalizations of the image conjecture and the Mathieu conjecture}, J. Pure and Appl.
Algebra 214 (7) (2010), 1200-1216.\\
\noindent [10] W. Zhao, {\em Mathieu Subspaces of Associative Algebras}, J. Algebra 350 (2012), 245-272.

\end{document}